\documentclass[10pt]{article}%

\usepackage{multirow}
\usepackage{tikz}
\usepackage{caption}
\usetikzlibrary{
  decorations.pathmorphing,%
  decorations.markings,%
  decorations.pathreplacing,%
  decorations.text,%
  calc,%
  patterns,%
  shapes.geometric,%
  arrows,
  decorations.shapes,
  positioning,
  plotmarks,
  shadings,
  math,
  intersections
}
\tikzset{>=latex} 
\tikzset{font=\small}
\tikzset{mark size=1.5pt, mark options=thin}
\tikzset{pin distance=4pt,
  every pin edge/.style={<-, thin, shorten <= -2pt}}

\definecolor{uipoppy}{RGB}{221,128,71}
\definecolor{uipaleblue2}{RGB}{179,196,215}
\definecolor{uiviolet}{RGB}{86,86,99}
\definecolor{uiblack}{RGB}{0, 0, 0}
\definecolor{azul}{RGB}{0,128,255}
\definecolor{verde}{RGB}{50,180,50}
\definecolor{uipaleblue}{RGB}{108,199,220}
\definecolor{light-gray}{gray}{0.8}
\definecolor{light-blue}{rgb}{0.53,.8,98}
\definecolor{green1}{RGB}{50,180,50}

\definecolor{jeffColor}{RGB}{102, 0, 204}
\definecolor{yaizaColor}{RGB}{0, 153, 153}
\definecolor{pale-verde}{RGB}{155,207,145}

\definecolor{periodColor}{RGB}{255, 167, 105}
\definecolor{dark-green}{RGB}{135, 194, 130}

\usepackage{bbm}
\usepackage{soul}

\usepackage{mathrsfs}
\DeclareMathAlphabet{\mathpzc}{OT1}{pzc}{m}{it}

\usepackage{color}
\usepackage{amsmath}
\usepackage{amsfonts}
\usepackage{amssymb}%
\usepackage{latexsym}

\newtheorem{theorem}{Theorem}[section]
\newtheorem{corollary}[theorem]{Corollary}
\newtheorem{definition}[theorem]{Definition}
\newenvironment{proof}[1][Proof]{\noindent \emph{#1.} }
{\hfill \ \rule{0.5em}{0.5em}}
\newtheorem{lemma}[theorem]{Lemma}
\newtheorem{proposition}[theorem]{Proposition}

\newtheorem{assumption}[theorem]{Assumption}

\numberwithin{equation}{section}
\numberwithin{table}{section}
\numberwithin{figure}{section}

\usepackage[a4paper]{geometry}
\newtheorem{remark}[theorem]{Remark}
\newtheorem{example}[theorem]{Example}
\newcommand{\noi}{\noindent}

\newcommand{\cA}{{\cal A}}

\newcommand{\cH}{{\cal H}}

\newcommand{\cO}{{\cal O}}

\newcommand{\bx}{\mathbf{x}}

\newcommand{\bxi}{\xi}
\newcommand{\bnu}{\boldsymbol{\nu}}

\newcommand{\ba}{\mathbf{a}}

\newcommand{\by}{\mathbf{y}}

    \newcommand\quotient[2]{
        \mathchoice
            {
                \text{\raise1ex\hbox{$#1$}\Big/\lower1ex\hbox{$#2$}}%
            }
            {
                #1\,/\,#2
            }
            {
                #1\,/\,#2
            }
            {
                #1\,/\,#2
            }
    }

\newcommand{\re}{{\rm e}}
\newcommand{\ri}{{\rm i}}
\newcommand{\rd}{{\rm d}}

\newcommand{\beq}{\begin{equation}}
\newcommand{\eeq}{\end{equation}}
\newcommand{\beqs}{\begin{equation*}}
\newcommand{\eeqs}{\end{equation*}}
\newcommand{\bit}{\begin{itemize}}
\newcommand{\eit}{\end{itemize}}
\newcommand{\ben}{\begin{enumerate}}
\newcommand{\een}{\end{enumerate}}
\newcommand{\bal}{\begin{align}}
\newcommand{\eal}{\end{align}}
\newcommand{\bals}{\begin{align*}}
\newcommand{\eals}{\end{align*}}
\newcommand{\bse}{\begin{subequations}}
\newcommand{\ese}{\end{subequations}}
\newcommand{\bpr}{\begin{proposition}}
\newcommand{\epr}{\end{proposition}}
\newcommand{\bre}{\begin{remark}}
\newcommand{\ere}{\end{remark}}
\newcommand{\bpf}{\begin{proof}}
\newcommand{\epf}{\end{proof}}
\newcommand{\ble}{\begin{lemma}}
\newcommand{\ele}{\end{lemma}}
\newcommand{\bco}{\begin{corollary}}
\newcommand{\eco}{\end{corollary}}
\newcommand{\bex}{\begin{example}}
\newcommand{\eex}{\end{example}}
\newcommand{\bth}{\begin{theorem}}
\newcommand{\enth}{\end{theorem}}

\newcommand{\Rea}{\mathbb{R}}
\newcommand{\Com}{\mathbb{C}}

\newcommand{\Oi}{{\Omega_-}}

\newcommand{\Oe}{{\Omega_+}}

\newcommand{\GR}{{\partial B_R}}

\newcommand{\eps}{\varepsilon}

\newcommand{\pdiff}[2]{\frac{\partial #1}{\partial #2}}

\newcommand{\gu}{\nabla u}

\newcommand{\vb}{\overline{v}}
\newcommand{\gvb}{\overline{\nabla v}}

\newcommand{\half}{\frac{1}{2}}

\newcommand{\tendi}{\rightarrow \infty}

\def\XXint#1#2#3{{\setbox0=\hbox{$#1{#2#3}{\int}$}
     \vcenter{\hbox{$#2#3$}}\kern-.5\wd0}}

\usepackage{hyperref}
\definecolor{myblue}{rgb}{0,0,0.6}
\hypersetup{colorlinks=true,linkcolor=myblue,citecolor=myblue,filecolor=myblue,urlcolor=myblue}
\newcommand*{\N}[1]{\left\|#1\right\|}

\allowdisplaybreaks[4]

\newcommand{\tfa}{\text{ for all }}
\newcommand{\tfor}{\text{ for }}

\newcommand{\tin}{\text{ in }}
\newcommand{\ton}{\text{ on }}
\newcommand{\tas}{\text{ as }}
\newcommand{\tand}{\text{ and }}
\newcommand{\tst}{\text{ such that }}

\newcommand{\vertiii}[1]{{\left\vert\kern-0.25ex\left\vert\kern-0.25ex\left\vert #1
    \right\vert\kern-0.25ex\right\vert\kern-0.25ex\right\vert}}

\newcommand{\DtN}{{\rm DtN}_k}

\definecolor{jwcol}{RGB}{27, 137, 18}  

\definecolor{dalcol}{rgb}{0.8,0,0}

\definecolor{escol}{rgb}{0,0,0.8}
\definecolor{jwcol}{rgb}{0,0.8,0}
\definecolor{estcol}{rgb}{0,0.5,0}
\definecolor{esnewcol}{rgb}{0,0.5,0}

\newcommand{\es}[1]{{\color{escol}{#1}}}

\newcommand{\Real}{\operatorname{Re}}
\newcommand{\Imag}{\operatorname{Im}}
\renewcommand{\Re}{\Real}
\renewcommand{\Im}{\Imag}

\newcommand{\supp}{{\rm supp}}
\newcommand{\comp}{{\rm comp}}

\newcommand{\RR}{\mathbb{R}}

\newcommand{\OR}{\Omega_R}

\newcommand{\mymatrix}[1]{\mathsf{#1}}
\newcommand{\MA}{{\mymatrix{A}}}
\newcommand{\MB}{{\mymatrix{B}}}
\newcommand{\MI}{{\mymatrix{I}}}
\newcommand{\Mzero}{{\mymatrix{0}}}
\newcommand{\SPD}{{\mathsf{SPD}}}

\newcommand{\Csol}{{C_{\rm sol}}}

\newcommand{\Ccont}{{C_{\rm cont}}}

\newcommand{\mythmname}[1]{\textbf{\emph{(#1)}}}

\newcommand{\pa}{\partial}

\newcommand{\domain}{\Omega}

\newcommand{\Id}{I}

\newcommand{\parao}{y}
\newcommand{\parat}{z}

\newcommand{\hatx}{\widehat{\bx}}
\newcommand{\GammaD}{\partial \Oi}
\newcommand{\CDtN}{{C_{\rm DtN}}}
\newcommand{\trace}{\gamma}
\newcommand{\tildeu}{\widetilde{u}}

\newcommand{\ind}{\mathbbm{1}}

\newcommand{\newy}{\by}

\begin{document}

\title{Wavenumber-explicit parametric holomorphy of Helmholtz solutions in the context of uncertainty quantification}

\author{E.~A.~Spence\footnotemark[1]\,\,, J.~Wunsch\footnotemark[2]}

\date{\today}

\footnotetext[1]{Department of Mathematical Sciences, University of Bath, Bath, BA2 7AY, UK, \tt E.A.Spence@bath.ac.uk }
\footnotetext[2]{Department of Mathematics, Northwestern University, 2033 Sheridan Road, Evanston IL 60208-2730, US, \tt jwunsch@math.northwestern.edu}

\maketitle


\begin{abstract}
A crucial role in the theory of uncertainty quantification (UQ) of PDEs is played by the regularity of the solution with respect to the stochastic parameters; indeed, a key property one seeks to establish is that the solution is holomorphic
with respect to (the complex extensions of) the parameters.
In the context of UQ for the high-frequency Helmholtz equation, a natural question is therefore:~how does this parametric holomorphy depend on the wavenumber $k$?

The recent paper \cite{GaKuSl:21} showed for a particular nontrapping variable-coefficient Helmholtz problem with affine dependence of the coefficients on the stochastic parameters that the solution operator can be analytically continued a distance $\sim k^{-1}$ into the complex plane.

In this paper, we generalise the result in \cite{GaKuSl:21} about $k$-explicit parametric holomorphy to a much wider class of Helmholtz problems with arbitrary (holomorphic) dependence on the stochastic parameters; we show that in all cases the region of parametric holomorphy decreases with $k$, and show how 
the rate of decrease with $k$ is 
dictated by whether the unperturbed Helmholtz problem is trapping or nontrapping.
We then give examples of both trapping and nontrapping problems where these bounds 
on the rate of decrease with $k$ of the region of parametric holomorphy are sharp, with the trapping examples coming from the recent results of \cite{GaMaSp:21}.

An immediate implication of these results is that the $k$-dependent restrictions imposed on the randomness in the analysis of quasi-Monte Carlo (QMC) methods in \cite{GaKuSl:21} arise from a genuine feature of the Helmholtz equation with $k$ large (and not, for example, a suboptimal bound).

\end{abstract}

\section{Introduction}

\subsection{Motivation:~wavenumber-explicit uncertainty quantification for the Helmholtz equation}

\paragraph{The importance of parametric analytic regularity in uncertainty quantification (UQ).}
The last $\sim$15 years has seen sustained interest in UQ of PDEs; i.e., the construction of algorithms (backed up by theory) for computing statistics of quantities of interest involving PDEs \emph{either} posed on a random domain \emph{or} having random coefficients.
These PDE problems can be posed in the abstract form 
\beq\label{eq:parametric}
P(\newy) u(\newy) = f
\eeq
where $P$ is a differential or integral operator, and $\newy$ is a vector of parameters governing the randomness. A crucial role in UQ theory is understanding regularity of $u$ with respect to the parameters $\newy$. Indeed, proving that $u$ is holomorphic with respect to (the complex extensions of) these parameters (see, e.g.,  \cite[Theorem 4.3]{CoDeSc:10}, \cite{CoDeSc:11},  \cite[Section 2.3]{KuSc:13}) is crucial for proving rates of convergence, independent of the number of the stochastic parameters,  
of 
 \bit
 \item stochastic collocation or sparse grid
schemes, see, e.g., \cite{ChCoSc:15,CaNoTe:16},
\item Smolyak quadratures, see, e.g., \cite{ZeDuSc:19, ZeSc:20}, 
\item quasi-Monte Carlo (QMC) methods, see, e.g., 
\cite{Sc:13, DiKuLeNuSc:14, DiKuLeSc:16, KuNu:16, HaPeSi:16}, and 
\item deep-neural-network approximations of the solution, see, e.g., \cite{ScZe:19, OpScZe:21, LoMiRuSc:21}.
\eit

\paragraph{UQ for the Helmholtz equation and $k$-explicit parametric regularity.}

Whilst a large amount of initial UQ theory concerned Poisson's equation $\nabla\cdot(\MA(\bx,\newy) \nabla u(\bx,\newy))= -f(\bx)$, there has been 
increasing interest in UQ of Helmholtz equation with (large) wavenumber $k$ 
(see, e.g.,
\cite{
FeLiLo:15, 
HiScScSc:15,
EsJe:20, PeSp:20, GrPeSp:21, 
BoNoPePr:20, GaKuSl:21})
and the time-harmonic Maxwell equations \cite{JeSc:17, JeScZe:17, FeLiLo:19, AlJeScZe:20}.
The Helmholtz equation with wavenumber $k$ and random coefficients is 
\beq\label{eq:Helmholtz}
k^{-2}
\nabla\cdot(\MA(\bx,\newy)\nabla u(\bx,\newy) ) +  n(\bx,\newy) u(\bx,\newy) = -f(\bx)
\eeq
where $\MA$ and $n$ depend on both the spatial variable $\bx$ and the stochastic variable $\newy$. 

Given the importance of holomorphy of $u$ with respect to $\newy$ for the parametric operator equation \eqref{eq:parametric}, 
a natural question in the context of the Helmholtz equation \eqref{eq:Helmholtz} is:
\begin{quotation}
\noi How does the holomorphy with respect to $\newy$ depend on $k$ as $k\to \infty$?
\end{quotation}
The recent paper \cite{GaKuSl:21} considers the Helmholtz equation
\eqref{eq:Helmholtz} posed for $\bx$ in a star-shaped domain $D$, with an impedance boundary condition on $\partial D$, and with 
\beq\label{eq:GKScoeff}
\MA\equiv \MI \quad\tand\quad n(\bx,\newy) = n_0(\bx)+\sum_{j=1}^\infty y_j \psi_j(\bx), 
\eeq
where $n_0$ and $\{\psi_j\}_{j=1}^\infty$ satisfy conditions so that, for every $\newy$, the coefficient $n$ does not trap geometric-optic rays and thus the solution operator has the best-possible dependence on $k$ (see \cite[Assumption A1]{GaKuSl:21}). The result \cite[Theorem 4.2]{GaKuSl:21} 
\footnote{Note that the $f$ in \cite[Theorem 4.2]{GaKuSl:21} is $k^2$ times the $f$ on the right-hand side of \eqref{eq:Helmholtz} because \cite{GaKuSl:21} 
considers the Helmholtz equation 
$\Delta u(\bx,\newy) +  k^2 n(\bx,\newy) u(\bx,\newy) = -f(\bx)$; see Remark \ref{rem:weighted} below for why we choose to write the Helmholtz equation as \eqref{eq:Helmholtz}).}
then shows that given $k_0>0$ there exist $C_0,C_j>0$ 
(with $C_j$ depending on $\|\psi_j\|_{W^{1,\infty}(D)}$) 
such that for all $k\geq k_0$, all $\newy$, and all finitely-supported multiindices $\alpha$,
\beq\label{eq:GaKuSl}
\N{\partial^\alpha_{\newy} u(\cdot,\newy)}_{L^2(D)} \leq C_0 
\bigg( \prod_j (k C_j)^{\alpha_j}\bigg)
|\alpha|!\, k\N{f}_{L^2(D)}.
\eeq
(In fact, \cite{GaKuSl:21} control a stronger norm of
$\partial^\alpha_{\newy} u(\cdot,\newy)$, and allow non-zero impedance data,
with the norm of this data then appearing with $\|f\|_{L^2(D)}$ on the right-hand side, but this is not important for our discussion here.)

If $C_j$ (which depends on $\|\psi_j\|_{W^{1,\infty}(D)}$) is independent of $k$ for all $j$, then
the bound \eqref{eq:GaKuSl} implies that, as a function of each $y_j$, the power series of $u$ has radius of convergence proportional to $k^{-1}$; this follows by bounding the Taylor series remainder using \eqref{eq:GaKuSl}.
Therefore, restricting attention to any finite set of the $\newy$ variables, $u$ (and hence also the solution operator $f\mapsto u$) has a holomorphic extension to a polydisc (i.e., a tensor product of discs in each $\newy$ coordinate) with radii $\sim k^{-1}$ and centred at the origin in the complex $\newy$ plane. 
Alternatively, to work with $\by$ variables whose size does not decrease with $k$, the analysis of QMC methods in \cite{GaKuSl:21} requires that $C_j \sim k^{-1}$ for all $j$, meaning that
$\|\psi_j\|_{W^{1,\infty}}$ decreases with $k$ (see \cite[Equations 5.6 and
5.7]{GaKuSl:21}). In either case, less random variation is
allowed as $k$ increases.

A similar result in the context of shape UQ for the Helmholtz equation is proved in \cite{HiScSp:22}. Indeed, 
for the Helmholtz transmission problem with parametric interface, \cite{HiScSp:22} proves that the solution operator is holomorphic in the interface parameters in a region $\sim k^{-1}$ when the basis functions describing the interface are independent of $k$.

\subsection{Informal summary of the results of this paper}

The main message of this paper is the following:
\begin{quotation}
\noi The region of parametric holomorphy of the Helmholtz solution operator \emph{decreases} as $k\tendi$, even when the solution operator has the best-possible dependence on $k$ (when the problem is nontrapping). When the problem is trapping (e.g., when the Helmholtz equation is posed outside an obstacle that traps geometric-optic rays), the region of parametric holomorphy can be exponentially-small for a sequence of $k$'s.
\end{quotation}
To show this, this paper contains the following results.
\ben
\item Lower bounds on the region of parametric holomorphy for a wide class of Helmholtz problems for arbitrary holomorphic perturbations, with the region of parametric holomorphy bounded in terms of the solution operator of the unperturbed problem;
see \S\ref{sec:set1} below. 
\item Two examples of the Helmholtz equation with a coefficient depending in an affine way on the parameter where, at least through an increasing sequence of $k$'s, the bounds in Point 1 are sharp. These two examples are the following:
\bit
\item A simple 1-d nontrapping example where explicit calculation shows that analytic extension to a ball of radius $k^{-1}$ is the largest possible; see \S\ref{sec:set2} below
\item The Helmholtz equation posed in the exterior of a strongly-trapping obstacle in $d\geq 2$, where analytic extension to a ball whose radius is exponentially-small in $k$ is the largest possible. This example follows directly from the recent results of \cite{GaMaSp:21}; see \S\ref{sec:set3} below.
\eit
\een 
The bounds in Point 1 generalise the parametric holomorphy result given by the bound \eqref{eq:GaKuSl} to a much wider class of scattering problems, and the 1-d nontrapping example in Point 2 implies sharpness of this parametric holomorphy result from \cite{GaKuSl:21}. 

The fact that the region of parametric holomorphy of the Helmholtz solution operator decreases with $k$ implies that 
when UQ algorithms relying on parametric holomorphy are applied to the Helmholtz equation
(at least without further modification) \emph{either} the performance  will degrade as $k\tendi$, \emph{or} constraints on the randomness that become more severe as $k\tendi$ must be imposed (as in \cite[Equations 5.6 and 5.7]{GaKuSl:21}) for the performance to be unaffected as $k\tendi$.

\subsection{$k$-explicit lower bounds on the region of parametric holomorphy for general Helmholtz problems}\label{sec:set1}

\paragraph{Informal description of the quantities in the statement of the result.} (For the precise definitions, see \S\ref{sec:definitions}.)
$\Oi$ is a bounded Lipschitz domain (allowed to be the empty set) such that its open complement
$\Oe:= \Rea^d\setminus \overline{\Oi}$ is connected. $\cA_0$ 
is the operator corresponding to the variational formulation of 
the Helmholtz equation
\beq\label{eq:Helmholtz1}
k^{-2}\nabla\cdot\big(\MA_0 \gu \big) + n_0 u = -f \quad \tin \domain_+, 
\eeq
with
\beq\label{eq:bc}
\text{ either } \quad \gamma u = 0 
\quad \text{ or }\quad \partial_{\nu, \MA_0} u =0
\ton \partial \Oi
\eeq
and satisfying the Sommerfeld radiation condition 
\beq\label{eq:src}
k^{-1}\pdiff{u}{r}(\bx) - \ri  u(\bx) = o \left( \frac{1}{r^{(d-1)/2}}\right)
\eeq
as $r:= |\bx|\tendi$, uniformly in $\hatx:= \bx/r$
(see Corollary \ref{cor:operators} below). In this definition, $\gamma$ is the trace operator and $\partial_{\nu, \MA} u$ is the conormal derivative of $u$ -- recall this is such that $\partial_{\nu, \MA}u = \nu \cdot\gamma (\MA\nabla u)$ when $u\in H^2$ near $\GammaD$, where $\nu$ is the outward-pointing unit normal vector on $\GammaD$.
Similarly, $\cA$ is the operator corresponding to the variational formulation of
\beqs
k^{-2}\nabla\cdot\big((\MA_0+ \MA_p) \gu \big) + (n_0+n_p) u = -f \quad \tin \domain_+,
\eeqs
\beqs
\text{ either } \quad \gamma u = 0 
\quad \text{ or }\quad \partial_{\nu, \MA_0+\MA_p} u =0
\ton \partial \Oi
\eeqs
with $u$ also satisfying the Sommerfeld radiation condition \eqref{eq:src} (i.e., we consider the Helmholtz equation \eqref{eq:Helmholtz} with $\MA= \MA_0+ \MA_p$ and $n=n_0+ n_p$).

We assume that both $\MA_0$ and $n_0$ are bounded above and below in $\Oe$, $\MA\equiv \MI$ outside a compact set contained in $B_R$ (the ball of radius $R$ centred at the origin), and $n\equiv 1$ outside $B_R$ (i.e., the support of $1-n$ can go up to $\partial B_R$, but the support of $\MI-\MA$ can't -- this is imposed so that, on $\partial B_R$, $\partial_{\nu, \MA} = \partial_\nu$) . 
Further conditions on $\MA_0$ are needed to ensure that the solution of \eqref{eq:Helmholtz1} exists and is unique; sufficient conditions are given in Theorem \ref{thm:Fred}, allowing discontinuous $\MA_0$. This set-up therefore covers scattering by either a Dirichlet or Neumann impenetrable obstacle and/or scattering by a penetrable obstacle (modelled by discontinuous $\MA_0$ and/or $n_0$).

We assume that both $\MA_p$ and $n_p$ are $L^\infty$ and supported in $B_R$, with the subscript ``$p$'' standing for ``perturbation''. We assume that $\MA_p$ and $n_p$ are holomorphic functions of a parameter $\newy$ in a subset of $\Com^N$, where $N$ is arbitrary; recall that, by Hartog's theorem on separate analyticity (see, e.g., \cite[Definition 2.1.1]{He:89}), this is equivalent to $\MA_p$ and $n_p$ being holomorphic functions of each $y_j$, $1\leq j\leq N$.

The operators $\cA_0$ and $\cA$ are defined using the
variational formulations of the problems on $\Omega_R:= \Oe \cap B_R$ (with the radiation condition realised via the exact Dirichlet-to-Neumann map on $\partial B_R$).
We work with the weighted norms
\beq\label{eq:weighted_norms}
\N{v}^2_{H^m_k(\domain_R)}:= \sum_{0\leq |\alpha|\leq m}k^{-2|\alpha|}\N{D^\alpha v}^2_{L^2(\domain_R)};
\eeq
the rationale for using these norms is that if a function $v$ oscillates with frequency $k$, then we expect $k^{-|\alpha|}|\partial^\alpha v|\sim  |v|$ (this is true, e.g., if $v(\bx) = \exp(\ri k \bx\cdot\ba)$ with $|\ba|=1$). 
Let $\cH$ be the Hilbert space defined by either 
\beq\label{eq:spaceEDP}
\cH:=  \big\{ v\in H^1(\domain_R) : \gamma v=0 \ton \partial \Oi\big\} \quad\text{ or }\quad \cH:= H^1(\domain_R),
\eeq
for, respectively, Dirichlet and Neumann boundary conditions, with the weighted norm 
\beq\label{eq:1knorm}
\N{v}^2_{\cH}:=\N{v}^2_{H^1_k(\domain_R)}:= k^{-2}\N{\nabla v}^2_{L^2(\domain_R)} + \N{v}^2_{L^2(\domain_R)}.
\eeq
Let $\cH^*$ be the space of antilinear functionals on $\cH$, with the norm
\beqs
\N{F}_{\cH^*}:= \sup_{v\in \cH\setminus\{0\}} \frac{|F(v)|}{\|v\|_{\cH}}.
\eeqs
The solution operator $\cA_0^{-1}$ is then well defined from $\cH^*\to
\cH$.

\begin{theorem}\mythmname{Lower bounds on regions of parametric holomorphy in terms of solution operator of unperturbed problem}\label{thm:2}
Suppose that $\Omega_-$, $\MA_0$, $n_0$, and $R_0$ satisfy Assumption \ref{ass:1} and $\MA_0$ is piecewise Lipschitz.
Suppose that 
$\MA_p$ and $n_p$ are holomorphic in $\newy$ 
for $\newy\in Y_0\subset \Com^N$ 
with values in $L^\infty(\OR, \Rea^d\times \Rea^d)$ and $L^\infty(\OR,\Rea)$, respectively.
Suppose further that $\MA_p(0)=\Mzero, n_p(0)=0$, and, for all $\newy \in Y_0$, 
$\supp \,n_p(\newy) \subset B_R$ and $\supp \,\MA_p(\newy) \Subset K \subset B_R$, where $K$ is independent of $\newy$.

(i)
Let $Y_{1}(k)\subset\Com^N$ be an open subset of $Y_0$ 
such that, for all $\newy\in Y_{1}(k)$, 
\beq\label{eq:condition1}
\N{\cA_0^{-1}(k)}_{\cH^* \to \cH}
\max\Big\{ \N{\MA_p(\newy)}_{L^\infty(\Omega_R)}, \N{n_p(\newy)}_{L^\infty(\Omega_R)}\Big\} \leq \frac{1}{2}.
\eeq
Then, for all $k>0$, the map $\newy\mapsto \cA^{-1}(k,\newy):\cH^*\to \cH$ is holomorphic for $\newy\in Y_{1}(k)$ with
\beq\label{eq:thm_bound1}
\N{\cA^{-1}(k,\newy)}_{\cH^* \to \cH} \leq 2 \N{\cA_0^{-1}(k)}_{\cH^* \to \cH}.
\eeq

(ii) Assume further that $\Omega_-$ is $C^{1,1}$, $\MA_0\in W^{1,\infty}(\OR, \SPD)$, and $\MA_p$ is holomorphic in $\newy$ 
for $\newy \in Y_0\subset \Com^N$
with values in $W^{1,\infty}(\OR,  \Rea^d\times \Rea^d)$.
Then given $k_0>0$ there exists $C>0$ (independent of $k$ and $\newy$) such that the following is true.
Let $Y_{2}(k)\subset\Com^N$ be an open subset of $Y_0$ such that, for all $\newy\in Y_{2}(k)$, 
\beq\label{eq:condition2}
\N{\cA_0^{-1}(k)}_{L^2(\OR) \to H^2_k(\OR)\cap \cH}
\max\Big\{C \N{\MA_p(\newy)}_{W^{1,\infty}(\Omega_R)}, \N{n_p(\newy)}_{L^\infty(\Omega_R)}\Big\} \leq \frac{1}{2}.
\eeq
Then, for all $k\geq k_0$, the map $\newy\mapsto 
\cA^{-1}(k,\newy):L^2(\OR) \to  H^2_k(\OR)\cap \cH$ is holomorphic for $\newy\in Y_{2}(k)$ with
\beqs
\N{\cA^{-1}(k,\newy)}_{L^2(\OR) \to  H^2_k(\OR)\cap \cH} \leq 2 \N{\cA_0^{-1}(k)}_{L^2(\OR) \to  H^2_k(\OR)\cap \cH}.
\eeqs
\end{theorem}

We make the following six remarks.

\bit
\item
Since $\MA_p(0)=\Mzero, n_p(0)=0$, the left-hand sides of \eqref{eq:condition1} and \eqref{eq:condition2} are zero when $y=0$, and thus the
sets $Y_{1}(k)$ and $Y_{2}(k)$ are non-empty. 
\item
Theorem \ref{thm:2} is a lower bound on the region where $\newy\mapsto \mathcal{A}^{-1}(k,\newy)$ is holomorphic since the theorem shows that (under the assumptions in Parts (i) and (ii), respectively) $Y_{1}(k)$ and $Y_{2}(k)$ lie inside this region, thus bounding this region from below.
\item
The matrix $L^\infty$ norm 
appearing in \eqref{eq:condition1} is defined by $\|\MB\|_{L^\infty(\OR)}:= {\rm ess sup}_{x\in \OR}\|\MB(\bx)\|_2$, where $\|\cdot\|_2$ denotes the spectral/operator norm on matrices induced by the Euclidean norm on vectors. The matrix $W^{1,\infty}$ norm 
appearing in \eqref{eq:condition2} is defined by $\|\MB\|_{L^\infty(\OR)} + \sum_{j=1}^d\|\partial_j \MB\|_{L^\infty(\OR)}$.
\item 
 In Part (ii), the assumption that $\MA_0\in W^{1,\infty}$ implies $H^2$ regularity of Helmholtz solutions with data in $L^2$; thus $\cA^{-1}_0(k): L^2(\OR)\to 
H^2_k(\OR)\cap\cH$ is well defined. The definitions of the weighted norms  implies that 
$\|\cA_0^{-1}(k)\|_{\cH^* \to \cH}$ has the same $k$-dependence as $\|\cA_0^{-1}(k)\|_{L^2(\OR) \to H^2_k(\OR)\cap \cH}$ (see \eqref{eq:H2bound} below).
\item
Holomorphy of $\cA^{-1}(k,\newy)$ for $\newy\in Y_{1}(k)$ or $Y_{2}(k)$ then implies derivative bounds similar to \eqref{eq:GaKuSl} (but with different constants and potentially-different $k$-dependence) by Cauchy's integral formula; see, e.g., \cite[Theorem 2.1.2]{He:89}, \cite[Proposition 2.2]{ZeSc:20}.
\item
We have restricted attention to $\newy\in \Com^N$, instead of considering, say, $\newy \in \Com^{\mathbb{N}}$, to avoid 
discussing the technicalities of what it means for a function of infinitely-many complex variables to be holomorphic (see, e.g., \cite{Mu:10}). 
We emphasise, however, that $N$ is arbitrary, and all the dependence of the conditions \eqref{eq:condition1} and \eqref{eq:condition2} on $N$ is contained in $\MA_p$ and $n_p$.
\eit

\paragraph{The $k$-dependence of Theorem \ref{thm:2}.}
The $k$-dependence of the conditions \eqref{eq:condition1} and \eqref{eq:condition2} is determined by the $k$-dependence of $\|\cA_0^{-1}(k)\|$. We now recap this $k$-dependence, omitting the spaces in this norm since each of 
\beqs
\|\cA_0^{-1}(k)\|_{\cH^* \to \cH},\,\, \|\cA_0^{-1}(k)\|_{L^2(\OR) \to H^2_k(\OR)\cap \cH}, \,\,\|\cA_0^{-1}(k)\|_{L^2(\OR)\to \cH}, \,\,\tand\,\, \|\cA_0^{-1}(k)\|_{L^2(\OR)\to L^2(\OR)}
\eeqs
has the same $k$-dependence thanks to the definition of the weighted norms \eqref{eq:weighted_norms} (see \eqref{eq:L2H11} and \eqref{eq:H2bound} below). In this next result we write $a\lesssim b$ if there exists $C>0$ (independent of $k$) such that $a\leq Cb$; below we write $a\sim b$ if $a\lesssim b$ and $b\lesssim a$.

\begin{theorem}[Informal statement of bounds on $\|\cA_0^{-1}(k)\|$ for large $k$]\label{thm:solutionoperator}

\

(i) $\|\cA_0^{-1}(k)\|\gtrsim k$.

(ii) $\|\cA_0^{-1}(k)\| \lesssim \exp(C k)$. 

(iii) If the problem is nontrapping, then $\|\cA_0^{-1}(k)\|\lesssim k$.

(iv) $\|\cA_0^{-1}(k)\|$ is polynomially bounded in $k$ for ``most'' frequencies; i.e., given $k_0,\delta,\eps>0$, there exists $J\subset [k_0,\infty)$ with $|J|\leq \delta$ such that
\beqs
\|\cA_0^{-1}(k)\|\lesssim k^{5d/2+2+\eps} \quad\tfa k\in [k_0,\infty)\setminus J.
\eeqs
\end{theorem}

Regarding (i):~when $\MA_0\equiv \MI$ and $n_0\equiv 1$ this lower bound can be proved by considering $u(\bx)= \re^{\ri k x_1} \chi(|\bx|/R)$ for $\chi \in C^\infty$ supported in $[0,1)$; see, e.g., \cite[Lemma 3.10]{ChMo:08}/\cite[Lemma 4.12]{Sp:14}.

Regarding (ii):~this exponential upper bound is proved in \cite[Theorem 2]{Bu:98} (for smooth $\MA_0$ and $n_0$, Dirichlet boundary conditions), \cite{Vo:00} 
(for smooth $\MA_0$ and $n_0$, Neumann boundary conditions), and \cite[Theorem 1.1]{Be:03} ($\MA_0$ and $n_0$ with one jump). This exponential bound is sharp through a sequence of $k$'s by \cite{Ra:71, St:00} (for certain smooth $n_0$), \cite{PoVo:99} (for certain $\MA_0$ and $n_0$ with a jump), and \cite[Proof of Theorem 2.8]{BeChGrLaLi:11} (for $\MA_0\equiv \MI, n_0\equiv 1$ and a certain Dirichlet $\Oi$).

Regarding (iii):~$\MA_0, n_0,$ and $\Omega_-$ are \emph{nontrapping} if the generalisations for variable $\MA_0$ and $n_0$ of geometric-optic rays propagating in a neighbourhood of $\Omega_-$ leave this neighbourhood in a uniform time. 
The concept of nontrapping is only rigorously well-defined when $\MA_0, n_0$ are both $C^{1,1}$ and $\Oi$ is $C^\infty$, and the bound $\|\cA^{-1}_0(k)\|\lesssim k$ is proved in this set-up in \cite{GaSpWu:20} (with the omitted constant given explicitly in terms of the longest ray in $\OR$). This bound for smooth $\MA_0,n_0,$ and $\Oi$ goes back to \cite{Va:75, MeSj:82}; see \cite[Theorem 4.4.3]{DyZw:19} and the discussion in \cite[\S4.7]{DyZw:19}.
The bound $\|\cA_0^{-1}(k)\|\lesssim k$ can also be established for certain discontinuous $\MA_0$ and $n_0$ which heuristically are ``nontrapping'', see \cite{CaPoVo:99, GrPeSp:19, MoSp:19}, although defining this concept rigorously for discontinuous coefficients is difficult.

Regarding (iv):~this was recently proved in \cite[Theorem 1.1]{LaSpWu:20}.

\begin{example}[Theorem \ref{thm:2} applied to the set-up in \cite{GaKuSl:21} of \eqref{eq:GKScoeff}]
Applied to the set-up in \cite{GaKuSl:21} of \eqref{eq:GKScoeff} (and assuming further that $\newy\in \Com^N$ for $N$ large), the condition \eqref{eq:condition1} is ensured if 
\beqs
\N{\cA_0^{-1}(k)}_{\cH^*\to \cH} \sum_{j=1}^N |y_j|\N{\psi_j}_{L^\infty} \leq \frac{1}{2}.
\eeqs
When $\MA_0, n_0,$ and $\Omega_-$ are nontrapping, Part (iii) of Theorem \ref{thm:solutionoperator} implies that, given $k_0>0$, there exists $\Csol>0$ such that 
$\|\cA_0^{-1}(k)\|_{\cH^*\to \cH}\leq \Csol k$ for all $k\geq k_0$ (with $\Csol$ given explicitly in terms of the longest ray in $\Omega_R$ by \cite{GaSpWu:20} when $\MA_0, n_0,$ and $\Omega_-$ are sufficiently smooth). Thus \eqref{eq:condition1} is ensured if 
\beq\label{eq:poly1}
k\sum_{j=1}^N |y_j|\N{\psi_j}_{L^\infty} \leq \frac{1}{2\Csol};
\eeq
i.e., applied to the set-up of \eqref{eq:GKScoeff}, Theorem \ref{thm:2} shows holomorphy of $\newy\mapsto \cA^{-1}(k,\newy):\cH^*\to \cH$ 
under the condition \eqref{eq:poly1}. 

There are (at least) two ways to proceed from here:~(i) assume that the norms of the basis functions $\{\psi_j\}_{j=1}^\infty$ decrease with $k$, and obtain holomorphy for $\by$ in a $k$-independent polydisc, or (ii) assume that the norms of $\{\psi_j\}_{j=1}^\infty$ are independent of $k$, and obtain holomorphy for $\by$ in a polydisc with radii $\sim k^{-1}$.

Regarding (i):~given $C>0$, if 
$k\sum_{j=1}^N\N{\psi_j}_{L^\infty} \leq (2\Csol C)^{-1}
$
then, for all $k\geq k_0$, $\newy\mapsto \cA^{-1}(k,\newy):\cH^*\to \cH$ is holomorphic in the polydisc $\{\by \in \Com^N : |y_j|\leq C \text{ for all } j\}$.
When  $y_j \in [-1/2,1/2]$, this condition ``$k\sum_{j=1}^N\|\psi_j\|_{L^\infty}$ sufficiently small'' was obtained in \cite[Appendix A, last line]{GaKuSl:21} (see the discussion below), and a similar condition imposed as  \cite[Equation 5.6]{GaKuSl:21}.

Regarding (ii):~given $C>0$, if $\sum_{j=1}^N\N{\psi_j}_{L^\infty} \leq (2\Csol C)^{-1}$,
then, for all $k\geq k_0$,
$\newy\mapsto \cA^{-1}(k,\newy):\cH^*\to \cH$ is holomorphic in the polydisc $\{\by \in \Com^N : |y_j|\leq C k^{-1}\text{ for all } j\}$.
\end{example}

\paragraph{The ideas behind Theorem \ref{thm:2} and relationship to \cite{Sc:13, GaKuSl:21}.}
The basic idea behind the proof of Theorem \ref{thm:2} is to treat $\cA$ as a perturbation of $\cA_0$ and use Neumann series.
This idea was also central to the abstract theory in \cite{Sc:13} of parametric operator equations with affine parameter dependence; this theory was then reviewed for the Helmholtz problem considered in \cite{GaKuSl:21} in \cite[Appendix A]{GaKuSl:21} (although the bound \eqref{eq:GaKuSl} is proved in \cite[\S4]{GaKuSl:21} by repeatedly differentiating the PDE with respect to $\newy$ and essentially applying the solution operator \footnote{A slight complication is that \cite{GaKuSl:21} use a non-standard variational formulation of the Helmholtz equation from \cite{GaMo:17a, GaMo:20, MoSp:14}.}). This type of perturbation argument is also implicitly used in the Helmholtz context in \cite{FeLiLo:15, GaGrSp:15, GrPeSp:21}.

The novelty of Theorem \ref{thm:2} is that it applies 
these abstract arguments to 
general Helmholtz problems (covering scattering by penetrable and impenetrable obstacles) with general holomorphic perturbations (in both the highest- and lowest-order term). Furthermore, Part (ii) of Theorem \ref{thm:2} gives conditions for parametric holomorphy of the solution operator mapping into $H^2$; recall that, first, $H^2$ spatial regularity is important 
for the analysis of finite-element approximations to solutions of the Helmholtz problem (see, e.g., \cite{LaSpWu:22} and the references therein, and, e.g., \cite[\S5]{Sc:13} for a discussion of this in a non-Helmholtz-specific UQ setting) and, second, 
the more regularity possessed by the image space of the solution operator as an analytic function, the better the result one can prove in the error analysis of multilevel QMC methods; see, e.g., \cite[\S4.3.1]{DiGaLeSc:17}.

\subsection{A 1-d nontrapping example showing $k$-explicit upper bounds on the region of parametric holomorphy through a sequence of $k$'s}
\label{sec:set2}

\paragraph{The 1-d Helmholtz problem.}

We consider the following 1-d Helmholtz operator on $\Rea^+$ with a zero Dirichlet boundary condition at $x=0$,
\beq\label{eq:Pdef}
P:=  k^{-2}\partial_x^2 +\left( 1 - \left(\frac12 + \parao \right) \ind_{[0,1]}(x)\right);
\eeq
i.e.,
\beq
\label{eq:PDE}
Pu(x):= 
\begin{cases}
\big( k^{-2}\partial_x^2 + 1\big)u(x), & x> 1,\\ 
\big( k^{-2}\partial_x^2 + (1/2- \parao) \big)u(x), & 0<x\leq 1.
\end{cases}
\eeq
The 1-d analogue of the Sommerfeld radiation condition \eqref{eq:src} is that 
\beqs
u(x) = \re^{\ri k x} \, \tfor x \gg 1  
\eeqs
and this is equivalent to the impedance boundary condition $(k^{-1}\partial_x - \ri)u|_{x=2}=0$ (i.e., in 1-d, the outgoing Dirichlet-to-Neumann map is an impedance boundary condition).

Observe that this problem falls into the framework used in Theorem \ref{thm:2} with $\Omega_- = \{0\}$, $\Omega_R= (0,2)$ (i.e., we truncate the ``exterior domain'' $(0,\infty)$ at $x=2$ and apply the outgoing Dirichlet-to-Neumann map there), $\MA_0\equiv \MI, \MA_p \equiv 0$, 
\beq\label{eq:n0p}
n_0(x):=
\begin{cases}
1 & x>1 \\
1/2 & 0\leq x\leq 1/2
\end{cases},
\quad \tand\quad n_p(x):= y\ind_{[0,1]}(x).
\eeq

\paragraph{Physical interpretation and link to \cite{GaKuSl:21}.}

The physical interpretation of $P$ is that there is a penetrable obstacle in $0<x\leq 1$, with 
$\parao$ controlling the wave speed inside the obstacle. We assume that $|\parao|\leq 1/4$ 
so that the wave speed inside the obstacle is bigger than the wave speed outside, but analogous results hold also in the opposite case.
This 1-d Helmholtz operator is therefore a simple model of 
the situation in \cite{GaKuSl:21} where the coefficient depends in an affine way on a parameter (here $\parao$).
Since the lower-order term in brackets in \eqref{eq:Pdef} is always $>0$, this 1-d problem is nontrapping; with the bound on the solution operator in Part (iii) of Theorem \ref{thm:solutionoperator} proved in, e.g., \cite[\S2.1.5]{Ch:15}, \cite[Theorem 1]{Ch:16}, \cite[Theorem 5.10]{GrSa:20}.

\paragraph{The solution operator and its meromorphic continuation.}
It is standard that the solution operator $\cA^{-1}(k,y): L^2(\OR)\to L^2(\OR)$ is a meromorphic family of operators
for $k\in \mathbb{C}$ with $\parao \in \RR$ fixed; see, e.g., \cite[\S\S2.2, 3.2, 4.2]{DyZw:19}.
The same arguments also show that, for fixed $k$, $\cA^{-1}(k,y):L^2(\OR)\to L^2(\OR)$ defines a
meromorphic family of operators as a function of
$\parao\in\mathbb{C}$; 
this is proved in a multi-dimensional setting in \cite[Lemma 1.12]{GaMaSp:21} (see \S\ref{sec:set3} below) and in the current 1-d setting in \S\ref{sec:simpleproofs}.

\begin{theorem}\textbf{\emph{(Non-zero solutions to $Pu=0$ with complex $\parao$ for certain increasing sequences of $k$'s)}}
\label{thm:1}
There exists $m_0\in \mathbb{Z}^+$ and $C_1,C_2>0$ such that if 
\beq\label{eq:k2}
k=2 \pi m \sqrt{2},
\eeq
with $m\geq m_0$ then there exists a non-zero solution $u\in H^1_{\rm loc}(\Rea)$ to  
\beqs
Pu = 0 \ton \Rea^+, \quad 
u(0)=0, \quad\tand \quad u(x) = \re^{\ri k x} \, \tfor x \geq 2,
\eeqs
with $\parao \in\mathbb{C}$ such that
\beq\label{eq:strip2}
C_1 \leq k |y| \leq C_2.
\eeq
\end{theorem}

\begin{corollary}\textbf{\emph{(Limits to analytic continuation of $\cA^{-1}(k,y)$ with respect to $\parao$ for certain increasing sequences of $k$s)}}\label{cor:1}
Under the assumptions of Theorem \ref{thm:1}, the solution operator $\cA^{-1}(k,y): L^2([0,2])\to L^2([0,2])$ cannot be analyticity continued as a function of $\parao$ to a ball centred at the origin of radius $\gg k^{-1}$ in the complex $\parao$-plane.
\end{corollary}

The proofs of  Theorem \ref{thm:1} and Corollary \ref{cor:1} are contained in \S\ref{sec:simpleproofs}.

\begin{corollary}[Theorem \ref{thm:2} applied to this 1-d example is sharp]
When applied to the Helmholtz problem with $P$ defined by \eqref{eq:Pdef}, the $k$-dependence of the 
condition \eqref{eq:condition1}
 in Theorem \ref{thm:2} is sharp through the sequence of $k$'s in Theorem \ref{thm:1}.
\end{corollary}

\bpf
As noted above, the Helmholtz problem with $P$ defined by \eqref{eq:Pdef} falls into the framework of Theorem \ref{thm:2} with $\MA_0\equiv \MI, \MA_p \equiv 0$, and $n_0$ and $n_p$ given by \eqref{eq:n0p}.
Since $\|n_p\|_{L^\infty} =|y|$ and $\|\cA_0^{-1}(k)\| \sim k$ by, e.g., \cite[\S2.1.5]{Ch:15}, \cite[Theorem 1]{Ch:16}, \cite[Theorem 5.10]{GrSa:20}, the bound \eqref{eq:condition1} implies that there exists $C>0$ (independent of $k$) such that the solution operator
$\cA^{-1}(k,y)$ is holomorphic in $y$ for $k|y|\leq C$. 
By Theorem \ref{thm:1}, when $k$ is given by \eqref{eq:k2}, there is a pole with $C_1 \leq k|y|\leq C_2$, with $C_1, C_2>0$ independent of $k$.
\epf

\bre[The idea behind Theorem \ref{thm:1}:~treating $y$ as a spectral parameter] \label{rem:idea_spectral}
A heuristic explanation of the limited analytic continuation in $y$ in Theorem \ref{thm:1} 
can be obtained by viewing the operator $P$ as a quantum Hamiltonian
$$
k^{-2} \pa_x^2 + E-V
$$ on $\Rea^+,$
where $k=\hbar^{-1}$ is the inverse Planck's constant, $E$ is a spectral parameter (with this notation indicating that it can also be considered as the energy of the system), and
$$
V= (1/2+y)\ind_{[0, 1]}(x)
$$
is a family of potentials, real-valued when $y \in \RR.$  It is well
known in such simple cases that energy \emph{reflects} off the potential
discontinuities at the boundary of the support, which yields a kind of
weak trapping.  When we then analytically continue the spectral
parameter $E$ across the continuous spectrum, this results in \emph{resonances},
i.e., poles of the meromorphic continuation of the solution operator
(see
\cite[Sections 2.2--2.3]{DyZw:19} for the description of this meromorphic
continuation).  The main observation used in the present paper is
that $y$ is behaving equivalently to the spectral parameter over a large
enough region, so we also get poles in $y.$

The resonant states, i.e., outgoing solutions to the ODE, can be
viewed as losing considerable energy from each reflection off the edge of
the barrier at $x=1.$ (Note that $E>V$, so there is no trapping of
orbits in the underlying classical system with Hamiltonian $p^2+V$ at
energy $E$.)  The energy then gets \emph{amplified} by
the nonzero imaginary part of $y$ as we propagate across $\supp(V),$
with the necessary balance of the effects of loss and amplification constraining the location of the
poles in $y.$
\ere

\bre[How would Theorem \ref{thm:1} change if the coefficient was continuous?]
\label{rem:conjecture}
The discontinuity of the coefficients in the example in Theorem \ref{thm:1} might make the reader wonder if this discontinuity is the sole
cause of the failure of holomorphic extension.  The analysis of
\cite{Be:82}, however, shows that \emph{some} reflection of energy still takes
place in the semiclassical ($k\to\infty$) limit if some derivative of
$V$ is discontinuous at the boundary of its support, or, indeed, even
from the necessary failure of analyticity of $V$ near the boundary of
its support (albeit more weakly). Similar reflection coefficients
arise in the analysis of scattering by $\delta$-potentials and of conic
diffraction, and this very weak trapping of energy turns out to yield resonances with imaginary part $\sim C k^{-1}
\log k$ \cite{GaSm:15, Ga:17, HiWu:20, DaMa:21}.   It thus seems reasonable to conjecture that the role of $y$
continues to be analogous to a spectral parameter in the setting of
\cite{Be:82} and that therefore the limits of analytic continuation in
$y$ should be no better than $C k^{-1} \log k$ even when $V$ is
$C^N$ -- see \cite{Re:58, Zw:87} for the analogous study of
resonance poles.
\ere

\subsection{A trapping example with region of parametric analyticity exponentially-small in $k$ (taken from \cite{GaMaSp:21})}\label{sec:set3}

We use the notation outlined in \S\ref{sec:set1} (and precisely defined in \S\ref{sec:definitions}).
Let 
\beq\label{eq:GMS1}
\MA_0\equiv \MI, \quad \MA_p \equiv 0, \quad n_0 \equiv 1, \,\,\tand\,\, n_p(\bx,y):= y \ind_{\OR}(\bx).
\eeq
We consider the case of zero Dirichlet boundary conditions on $\partial \Omega_-$ (i.e., the first boundary condition in \eqref{eq:bc}). 

\begin{lemma}\mythmname{Meromorphy of the solution operator as a function of $y$ \cite[Lemma 1.12]{GaMaSp:21}}\label{lem:meromorphy}
Let $\MA$ and $n$ be defined by \eqref{eq:GMS1} and let the space $\cH$ be
defined by the first equation in \eqref{eq:spaceEDP}. Then, for all
$k>0$, $y\mapsto \cA^{-1}(k,y): L^2(\OR)\to L^2(\OR)$ is a meromorphic
family of operators for $y\in \Com$.
\end{lemma}

It is well-known (see Lemma \ref{lem:specific} below) that domains with strong trapping have so-called \emph{quasimodes}, and the results of \cite{GaMaSp:21} are most-easily stated using this concept.

\begin{definition}[Quasimodes]\label{def:quasimodesh}
A \emph{family of quasimodes of quality $\eps(k)$}
is a sequence $\{(u_\ell,k_\ell)\}_{\ell=1}^\infty\subset H^2(\Omega_R)\cap \cH
\times \mathbb{R}$ such that  $k_\ell\to\infty$ as $\ell \tendi$ and there is a compact subset $\mathcal{K}\Subset B_R$ such that, for all $\ell$, $\supp\, u_\ell \subset \mathcal{K}$,
\beqs
\N{(-k^{-2}\Delta -1) u_\ell}_{L^2(\Omega_R)} \leq \eps(k_\ell), \quad\tand\quad\N{u_\ell}_{L^2(\Omega_R)}=1.
\eeqs
\end{definition}

\begin{theorem}\mythmname{From quasimodes to poles of the solution operator in $y$ \cite[Theorem 2.2]{GaMaSp:21}}
\label{thm:GMS}
Let $\alpha> 3(d+1)/2$.
Suppose there exists a family of quasimodes {in the sense of Definition \ref{def:quasimodesh}} 
such that the quality $\eps(k)$ satisfies
\beq\label{eq:epslowerbound}
\eps(k) =o( k^{-1-\alpha}) \quad\tas k \to\infty.
\eeq
Then there exists $k_0>0$ (depending on $\alpha$) such that, if $\ell$ is such that $k_\ell\geq k_0$ then there exists $y_\ell\in\mathbb{C}$ with 
\beq\label{eq:pole}
|{y_\ell}| \leq k_\ell^{\alpha} \eps(k_\ell)
\eeq
such that $y\mapsto \cA^{-1}(k,y):{L^2(\OR)\to L^2(\OR)}$ has a pole at $y_\ell$.
\end{theorem}

In \S\ref{sec:GMS} we explain how Lemma \ref{lem:meromorphy} and Theorem \ref{thm:GMS} follow from the results of \cite{GaMaSp:21}.
We note that \cite{GaMaSp:21} in fact proves the stronger result that quasimodes imply poles in $y$ with the same multiplicities; see \cite[Theorems 1.8 and 2.4]{GaMaSp:21}.

The following lemma gives three specific cases when the assumptions of Theorem \ref{thm:GMS} hold; this result uses the notation that 
$B = \cO(k^{-\infty})$ as $k\tendi$ if, given $N>0$, there exists $C_N$ and $k_0$ such that $|B|\leq C_N k^{-N}$ for all $k\geq k_0$, i.e.~$B$ decreases superalgebraically in $k$.

\begin{lemma}[Specific cases when the assumptions of  Theorem \ref{thm:GMS} hold]\label{lem:specific}

\

(i) Let $d=2$. Given $a_1>a_2>0$, let 
\beq\label{eq:ellipse}
E:= \left\{(x_1,x_2) \, : \, \left(\frac{x_1}{a_1}\right)^2+\left(\frac{x_2}{a_2}\right)^2<1\right\}.
\eeq
If $\partial \Omega_-$ coincides with the boundary of $E$ in the neighborhoods of the points \((0,\pm a_2)\), 
and if $\Omega_+$ contains the convex hull of these neighbourhoods, 
then the assumptions of Theorem \ref{thm:GMS} hold with 
\beqs
\eps(k)= \exp( - C_1 k)
\eeqs
for some $C_1>0$ (independent of $k$).\footnote{
In \cite[Theorem 2.8]{BeChGrLaLi:11}, $\Omega_+$ is assumed to contain the whole ellipse $E$.
However, inspecting the proof, we see that the result remains unchanged if $E$ is replaced with the convex hull of the neighbourhoods of 
$(0,\pm a_2)$. Indeed, the idea of the proof is to consider a family of eigenfunctions of the ellipse localising around the periodic orbit $\{(0,x_2) : |x_2|\leq a_2\}$.
}

(ii) Suppose $d\geq 2$, $\Gamma_D\in C^\infty$, and $\Omega_+$
contains an elliptic closed orbit (roughly
  speaking, a trapped ray that is stable under
  perturbation\footnote{More precisely, 
  the eigenvalues of the linearized
  Poincar\'e map have moduli $\leq 1$.}) such
that (a) $\partial\Omega_-$ is analytic in a neighbourhood of the ray
and (b) the ray satisfies the stability and nondegeneracy condition
\cite[(H1)]{CaPo:02}. If $q >11/2$ when $d=2$ and $q>2d+1$ when
$d\geq 3$, then the assumptions of Theorem \ref{thm:GMS} hold with
\beqs \eps(k)= \exp( - C_2 k^{1/q}) \eeqs for some $C_2>0$
(independent of $k$).

(iii) Suppose there exists a sequence of resonances $\{\lambda_\ell\}_{\ell=1}^\infty$ of the exterior Dirichlet problem with
\beqs
0\leq -\Im \lambda_\ell = \mathcal{O}\big(|\lambda_\ell|^{-\infty}\big)  \quad\tand \quad \Re \lambda_\ell \tendi \quad\tas\quad \ell \tendi.
\eeqs
Then there exists a family of quasimodes of quality 
\beqs
\eps(k)=\mathcal{O}(k^{-\infty}),
\eeqs
and thus the assumptions of Theorem \ref{thm:GMS} hold.
\end{lemma}

\bpf[References for the proof]
Parts (i) and (ii) are via the quasimode constructions of 
\cite[Theorem 2.8, Equations 2.20 and 2.21]{BeChGrLaLi:11} and \cite[Theorem 1]{CaPo:02} for obstacles whose exteriors support elliptic-trapped rays. Part (iii) is via the ``resonances to quasimodes'' result of \cite[Theorem 1]{St:00}; recall that the resonances of the exterior Dirichlet problem are the poles of the meromorphic continuation of the solution operator 
from $\Im k\geq 0$ to $\Im k<0$; see, e.g., \cite[Theorem 4.4 and Definition 4.6]{DyZw:19}.
\epf

\begin{figure}[h!]
    \centering
    \includegraphics[width=0.6\textwidth]{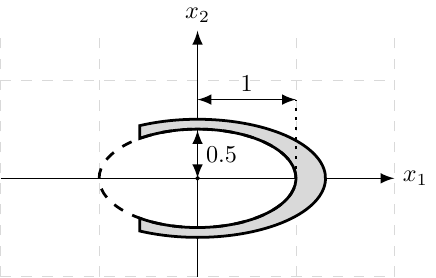}
    \caption{An example of an $\Omega_-$ (shaded) satisfying Part (i) of Lemma \ref{lem:specific} with $a_1=1$ and $a_2=1/2$. The part of the boundary of the ellipse \eqref{eq:ellipse} that is not part of $\partial \Omega_-$ is denoted by a dashed line}\label{fig:ellipse}
   \end{figure}

\begin{corollary}\mythmname{Theorem \ref{thm:1} applied to certain trapping $\Oi$} 
\label{cor:sharp2}
Suppose that $\Oi$ is as in Part (i) of Lemma \ref{lem:specific}, with $\partial \Oi$ additionally $C^\infty$ (i.e., $\Oi$ can be the obstacle in Figure \ref{fig:ellipse} with the corners smoothed).
Let $\MA$ and $n$ defined by \eqref{eq:GMS1} and let the space $\cH$ be defined by the first equation in \eqref{eq:spaceEDP}.
Let $(k_\ell)_{\ell=0}^\infty$ be the wavenumbers in the quasimode for this $\Oi$ (which exists by Lemma \ref{lem:specific}).

Then there exists $C_3>0$ and $\ell_0$ such that, for $\ell\geq \ell_0$, the map $y\mapsto \cA_0^{-1}(k_\ell,y)$ has a pole at $y_\ell$ with
\beq\label{eq:pole2}
|y_\ell | \leq \exp(-C_3 k_\ell).
\eeq
\end{corollary}

Corollary \ref{cor:sharp2} shows that Theorem \ref{thm:2} applied to the set-up of Corollary \ref{cor:sharp2} is sharp through the sequence $k_\ell$ (for $\ell$ sufficiently large). Indeed, in the set-up of Corollary \ref{cor:sharp2}, $\|n_p\|_{L^\infty}=|y|$, and $\|\cA_0^{-1}\|\lesssim \exp (C_4 k)$ for some $C_4>0$ (independent of $k$) by the result of \cite{Bu:98} in Part (ii) of Theorem \ref{thm:solutionoperator}. The condition \eqref{eq:condition1} for holomorphy therefore becomes 
\beqs
|y|\lesssim \exp(-C_4 k),
\eeqs
which is then sharp in its $k$-dependence when $k=k_\ell$ by the presence of the pole satisfying \eqref{eq:pole2}.

\paragraph{Outline of the rest of the paper.}
\S\ref{sec:definitions} defines the Helmholtz sesquilinear form and variational problem. 
\S\ref{sec:proof2} proves 
Proof of Theorem \ref{thm:2}.
\S\ref{sec:simpleproofs} proves Theorem \ref{thm:1} and Corollary \ref{cor:1}.

\section{Definition of the Helmholtz sesquilinear form and variational problem}\label{sec:definitions}

\subsection{Notation and assumptions on the domain and the coefficients}

\textbf{Notation:} $L^p(\domain)$ denotes complex-valued $L^p$ functions on a Lipschitz open set $\domain$.
When the range of the functions is not $\Com$, it will be given in the second argument; e.g. 
$L^\infty(\domain,\Rea^{d\times d})$ denotes the space of $d\times d$ matrices with each entry a real-valued $L^\infty$ function on $\domain$.
We use $\gamma$ to denote the trace operator $H^1(\domain)\rightarrow H^{1/2}(\partial \domain)$ and $\partial_\nu$ to denote the normal derivative trace operator $H^1(\domain,\Delta)\rightarrow H^{-1/2}(\partial \domain)$, where $H^1(\domain,\Delta):= \{ v \in H^1(\domain) : \Delta v\in L^2(\Omega)\}$.

\begin{assumption}
\mythmname{Assumptions on the domain and coefficients}
\label{ass:1}

(i) $\Omega_- \subset\Rea^d, d=2,3,$ is a bounded open Lipschitz set such that its open complement $\Omega_+:= \Rea^d\setminus \overline{\Omega_-}$ is connected. 

(ii) $\MA_0 \in L^\infty(\Omega_+ , \SPD)$
(where $\SPD$ is the set of $d\times d$ real, symmetric, positive-definite matrices)
is  such that $\supp(\MI- \MA_0)$ is bounded and
there exist $0<A_{0,\min}\leq A_{0,\max}<\infty$ such that, 
for all $\bxi\in \Rea^d$,
\beqs
 A_{0,\min}|\bxi|^2 \leq 
\big(\MA_0(\bx)\bxi\big)\cdot \bxi
 \leq A_{0,\max}|\bxi|^2 \quad\text{ for every }\bx \in \Omega_+.
\eeqs

(iii) $n_0\in L^\infty(\Omega_+,\Rea)$ is such that $\supp(1-n_0)$ is bounded and there exist $0<n_{0,\min}\leq n_{0,\max}<\infty$ such that
\beqs
n_{0,\min} \leq n_0(\bx)\leq n_{0,\max}\quad \text{ for almost every } \bx \in \Omega_+.
\eeqs

(iv) $R>0$ is such that $\overline{\Omega_-}\cup \supp(\MI- \MA)\Subset B_R$ and $\supp(1-n) \subset B_R$, where $B_R$ denotes the ball of radius $R$ about the origin.
\end{assumption}

 Let $\Omega_R:= \Omega_+\cap B_R$, and let $\Gamma_R:= \partial B_R$.

\subsection{The Sommerfeld radiation condition and the Dirichlet-to-Neumann map}\label{sec:src}

We say that $u\in C^1(\Rea^d\setminus B_{R'})$ for some $R'>0$ is \emph{outgoing} if it satisfies the Sommerfeld radiation condition \eqref{eq:src}.

\ble\mythmname{Explicit Helmholtz solution in the exterior of a ball}
Given $g\in H^{1/2}(\Gamma_R)$, the outgoing solution $v$ of 
\beq\label{eq:DtN1}
(-k^{-2}\Delta - 1)u=0 \quad\tin \Rea^d \setminus \overline{B_R} \quad\tand\quad \gamma u =g \ton \Gamma_R
\eeq
is unique and is given when $d=2$ by 
\beq\label{eq:sov}
u(r,\theta) = \frac{1}{2\pi} \sum_{n=-\infty}^\infty \frac{H^{(1)}_n(kr)}{ H^{(1)}_n(kR)} \exp(\ri n \theta) \widehat{g}(n),
\quad\text{ where }\quad
\widehat{g}(n):= \int_0^{2\pi} \exp(-\ri n\theta) g(R, \theta)\, \rd \theta
\eeq
(an analogous expression is available for $d=3$ -- see, e.g., 
\cite[Theorem 2.37]{KiHe:15}, \cite[\S3]{ChMo:08}, \cite[\S3]{MeSa:10}).
\ele

\bpf[References for the proof]
The uniqueness result is Rellich's uniqueness theorem; see, e.g., \cite[Theorem 3.13]{CoKr:83}. 
For the proof that \eqref{eq:sov} is an outgoing solution to \eqref{eq:DtN1}, see, e.g., \cite[Theorem 2.37]{KiHe:15}.
\epf

\

Given $g\in H^{1/2}(\Gamma_R)$, let $v$ be the outgoing solution of \eqref{eq:DtN1}.
Define the map $\DtN: H^{1/2}(\GR)\to H^{-1/2}(\GR)$ by 
\beq\label{eq:DtNdef}
\DtN g := k^{-1} \partial_{\bnu} v, 
\eeq
where $\bnu:= \bx/R= \hatx$ (i.e., $\bnu$ is the outward-pointing unit normal vector to $B_R$),
so that, when $d=2$, by \eqref{eq:sov},
\beq\label{eq:DtN}
\DtN g(\theta) = \frac{1}{2\pi} \sum_{n=-\infty}^\infty \frac{H^{(1)'}_n(kR)}{ H^{(1)}_n(kR)} \exp(\ri n \theta) \widehat{g}(n).
\eeq

\ble\mythmname{Key properties of $\DtN$}\label{lem:DtN}

(i) Given $k_0,R_0>0$ there exists $\CDtN_1= \CDtN_1(k_0 R_0)$ such that for all $k\geq k_0$ and $R\geq R_0$,
\beqs
\big|\big\langle \DtN\trace u, \trace v\rangle_{\partial B_R}\big\rangle\big| \leq k\,\CDtN_1 \N{u}_{H^1_k(\OR)}  \N{v}_{H^1_k(\OR)} \quad\tfa u,v \in H^1(\OR).
\eeqs

(ii) 
\beq\label{eq:DtN_imag}
\Im \big\langle \DtN \phi,\phi\big\rangle_{\GR} >0 \quad\tfa \phi \in H^{1/2}(\GR)\setminus \{0\}.
\eeq

(iii) 
\beq\label{eq:CDtN2}
- \Re \big\langle \DtN \phi,\phi\big\rangle_{\GR} \geq 0
 \quad\tfa \phi\in H^{1/2}(\GR).
\eeq
\ele

\bpf[References for the proof of Lemma \ref{lem:DtN}]
Each of (i), (ii), and (iii) are proved using the expression for $\DtN$ in terms of Bessel and Hankel functions (i.e., \eqref{eq:DtN} when $d=2$) in \cite[Lemma 3.3]{MeSa:10} (see also \cite[Theorem 2.6.4]{Ne:01} and \cite[Lemma 2.1]{ChMo:08} for \eqref{eq:DtN_imag} and \eqref{eq:CDtN2}).
\epf

\subsection{The Helmholtz sesquilinear form and associated operators}

\begin{definition}\mythmname{Helmholtz sesquilinear forms and associated operators}
Given $\MA_0,n_0, \Oi,$ and $R$ satisfying Assumption \ref{ass:1}, let 
\beq\label{eq:EDPa}
a_0(u,v):= \int_{\domain_R} 
\Big(k^{-2}(\MA_0 \nabla u)\cdot\gvb
 - n_0u \vb\Big) - k^{-1}\big\langle \DtN \gamma u,\gamma v\big\rangle_{\Gamma_R}.
\eeq
Given $\MA_p \in L^\infty(\domain_R, \Rea^d\times \Rea^d)$ and $n_p \in L^\infty(\domain_R,\Rea)$ with $\supp \,\MA_p \cup \supp \,n_p\Subset B_R$, let 
\beq\label{eq:a_p}
a_p(u,v):= \int_{\domain_R} 
\Big(k^{-2}(\MA_p\nabla u)\cdot\gvb
 - n_p u \vb\Big)
\eeq
(the subscript ``p'' standing for ``perturbation").
Let 
\beqs
a:= a_0 + a_p.
\eeqs
\end{definition}

The following lemma is proved using Parts (i) and (ii) of Lemma \ref{lem:DtN}, the Cauchy-Schwarz inequality, and  the definition of $\|\cdot\|_{H^1_k(\OR)}$ \eqref{eq:1knorm}.

\ble\mythmname{Properties of $a_0$ and $a_p$}\label{lem:propa}

(i) (Continuity of $a_0$) 
Given $k_0, R_0>0$ there exists $\Ccont>0$ such that for all $k\geq k_0$ and $R\geq R_0$, 
\beqs
|a_0(u,v)|\leq \Ccont \N{u}_{H^1_k(\Omega_R)} \N{v}_{H^1_k(\Omega_R)} \quad\tfa u, v \in H^1(\Omega_R).
\eeqs

(ii) (Continuity of $a_p$) 
\beqs
|a_p(u,v)|\leq \max\big\{ \N{\MA_p}_{L^\infty(\Omega_R)}, \N{n_p}_{L^\infty(\Omega_R)}\big\} \N{u}_{H^1_k(\OR)} \N{v}_{H^1_k(\OR)}
 \quad\tfa u, v \in H^1(\Omega_R).
\eeqs

(ii) (G\aa rding inequality for $a_0$) 
\begin{align*}
\Re a_0(v,v) &\geq  A_{0,\min} \N{v}^2_{H^1_k(\domain_R)} - \big(n_{0,\max}+ A_{0,\min}\big)\N{v}^2_{L^2(\domain_R)} \quad \tfa v\in H^1(\Omega_R).
\end{align*}
\ele

Recall the definition of the space $\cH$ \eqref{eq:spaceEDP} and $\|\cdot\|_{H^1_k(\OR)}$ \eqref{eq:1knorm}. 
Lemma \ref{lem:propa} combined with, e.g., \cite[Lemma 2.1.38]{SaSc:11} implies the following corollary.

\begin{corollary}\mythmname{The operators $\cA_0$ and $\cA_p$}\label{cor:operators}
There exist unique $\cA_0, \cA_p : \cH \to \cH^*$ such that, with $j$ equal either $0$ or $p$,
\beq\label{eq:Ajdual}
\langle \cA_j u, v\rangle_{\cH^* \times \cH} = a_j(u,v) \quad\tfa u,v \in \cH,
\eeq
\beq\label{eq:normboundA_p}
\N{\cA_0}_{\cH \to\cH^*}\leq \Ccont \quad\tand\quad \N{\cA_p}_{\cH \to\cH^*}\leq \max\Big\{ \N{\MA_p
}_{L^\infty(\Omega_R)}, \N{n_p
}_{L^\infty(\Omega_R)}\Big\}.
\eeq
\end{corollary}

\bre\mythmname{Approximating $\DtN$}
Implementing the operator $\DtN$ appearing in $a(\cdot,\cdot)$ \eqref{eq:EDPa} is computationally expensive, and so in practice one seeks to approximate this operator by \emph{either} imposing an absorbing boundary condition on $\GR$, \emph{or} using a perfectly-matched layer (PML), \emph{or} using boundary integral equations (so-called ``FEM-BEM coupling''). 
Recent $k$-explicit results on the error incurred (on the PDE level) by approximating $\DtN$ by absorbing boundary conditions or PML can be found in \cite{GaLaSp:21} and \cite{GaLaSp:21a}, respectively.
\ere

\subsection{The variational formulation and the solution operator}

\begin{definition}\mythmname{Variational formulation}\label{def:EDPvar}
Given $\Omega_-$, $\MA_0$, $n_0$, and $R_0$ satisfying Assumption \ref{ass:1} and 
$F\in \cH^*$, 
\beq\label{eq:EDPvar}
\text{ find } \tildeu \in \cH \tst \quad a_0(\tildeu,v)=F(v) \quad \tfa v\in \cH
\eeq
(i.e., $\cA_0 \tildeu =F$ in $\cH^*$).
\end{definition}

\ble\mythmname{From the variational formulation to the standard weak form}
Suppose that $\Omega_-$, $\MA_0$, $n_0$, and $R_0$ satisfy Assumption \ref{ass:1} and, in addition, 
$\MA_0 \in C^{0,1}(\OR, \SPD)$. Given $f \in L^2(\OR)$ with $\supp \, f \subset B_R$, let $F(v):=\int_{\OR} f\,\overline{v}$.

If $u\in H^1_{\rm{loc}}(\domain_+)$ is an outgoing solution of \eqref{eq:Helmholtz1} and \eqref{eq:bc} 
(where the PDE is understood in the standard weak sense), then $u|_{B_R}$ is a solution of the variational problem \eqref{eq:EDPvar}. 

Conversely, if $\widetilde{u}$ is a solution of this variational problem, then there exists a solution $u$ of \eqref{eq:Helmholtz}-\eqref{eq:src} such that $u|_{B_R}= \widetilde{u}$.
\ele

\bpf[Sketch of the proof]
This follows from Green's identity (see, e.g., \cite[Lemma 4.3]{Mc:00}) and the definitions of $\DtN$ and $a_0(\cdot,\cdot)$ \eqref{eq:EDPa}.
Note that it is crucial that $\MA\equiv \MI$ in a neighbourhood of $\GR$, so that the conormal derivative $\partial_{\nu, \MA}$ equals $\partial_{\nu}$ on $\GR$ 
(with the latter appearing in the definition of $\DtN$ \eqref{eq:DtNdef}).
\epf

\begin{theorem}[Invertibility of $\cA_0$]\label{thm:Fred}
Suppose that $\Omega_-$, $\MA_0$, $n_0$, and $R_0$ satisfy Assumption \ref{ass:1} and, in addition, $\MA_0$ is piecewise Lipschitz
(in the sense described in \cite[Proposition 2.13]{LiRoXi:19}).
Then $\cA_0^{-1} : \cH^* \to \cH$ is bounded.
\end{theorem}

\bpf[Sketch of the proof]
Under these assumptions, the Helmholtz equation satisfies a unique continuation principle (UCP); indeed, a UCP for $n_0\in L^{3/2}$ is proved in  \cite{JeKe:85,  Wo:92}, and a UCP for piecewise Lipschitz $\MA_0$ is proved in \cite[Proposition 2.13]{LiRoXi:19} using the results of \cite{BaCaTs:12}.
The UCP combined with Part (ii) of Lemma \ref{lem:DtN} imply that the solution of the variational problem \eqref{eq:EDPvar} is unique; i.e. $\cA_0$ is injective. 

Since the sesquilinear form is continuous (by Lemma \ref{lem:propa}) and satisfies a G\aa rding inequality, Fredholm theory implies 
that existence of a solution to the variational problem and continuous dependence of the solution on the data both follow from uniqueness; see, e.g., \cite[Theorem 2.34]{Mc:00}, \cite[\S6.2.8]{Ev:98}. 
\epf

\ble\mythmname{Norms of solution operators between different spaces}\label{lem:L2H1}
Suppose that $\Omega_-$, $\MA_0$, $n_0$, and $R$ satisfy Assumption \ref{ass:1}. Then
\beq\label{eq:L2H11}
\N{\cA_0^{-1}(k)}_{\cH^*\to \cH} \leq \frac{\big(1+ 2 n_{0,\max} \N{\cA_0^{-1}(k)}_{L^2(\OR)\to \cH}\big)}{\min\{A_{0,\min}, n_{0,\min}\}},
\eeq
and 
\beq\label{eq:L2H12}
\N{\cA_0^{-1}(k)}_{L^2(\OR)\to \cH} \leq \sqrt{\frac{3 n_{0,\max}}{2 A_{0,\min}} +1 }
\N{\cA_0^{-1}(k)}_{L^2(\OR)\to L^2(\OR)} + \frac{1}{2 n_{0,\min} k^2}.
\eeq
\ele

\bpf[References for the proof]
For \eqref{eq:L2H11}, see, e.g., \cite[Text between Lemmas 3.3 and 3.4]{ChMo:08} or \cite[Lemma 5.1]{GrPeSp:19}.
For \eqref{eq:L2H12}, see, e.g., \cite[Lemma 3.10, Part (i)]{GrPeSp:19}.
\epf

\begin{theorem}\mythmname{$H^2$ regularity}\label{thm:H2}
Suppose that $\Omega_-$, $\MA_0$, $n_0$, and $R_0$ satisfy Assumption \ref{ass:1} and, in addition, $\Omega_-$ is $C^{1,1}$ and $\MA_0$ is $W^{1,\infty}$. 
Then $\cA_0^{-1} : L^2(\OR) \to H^2(\OR)\cap \cH$ is bounded, and there exists $C>0$ (independent of $k$ but depending on $\MA_0$ and $n_0$) such that
\beq\label{eq:H2bound}
\N{\cA_0^{-1}(k)}_{L^2(\OR) \to H^2_k(\OR)\cap \cH} \leq C \N{\cA_0^{-1}(k)}_{\cH^*(R+1)\to \cH(R+1)}
\eeq
where
(i) on the space $H^2_k(\OR)\cap \cH$
 we use the $H^2_k(\OR)$ norm defined by \eqref{eq:weighted_norms}, and (ii) we write $\cH^*(R+1)$ for the space $\cH$ defined by \eqref{eq:spaceEDP} with $R$ replaced by $R+1$. 
\end{theorem}

\bpf Given $f\in L^2(\OR)$, let $u$ be the outgoing solution of \eqref{eq:Helmholtz1}.
By elliptic regularity (see, e.g., \cite[Theorem 4.18]{Mc:00}, \cite[\S6.3.2]{Ev:98}), there exists $C>0$ (depending on the $W^{1,\infty}$ norm of $\MA_0$) such that 
\beqs
|u|_{H^2(\OR)} \leq C \Big( \N{\nabla \cdot (\MA_{0}\nabla u)}_{L^2(\Omega_{R+1})}+ R^{-1}\N{\nabla u}_{L^2(\Omega_{R+1})}+ R^{-2}\N{u}_{L^2(\Omega_{R+1})}\Big).
\eeqs
The result then follows by bounding the right-hand side in terms of $\|f\|_{L^2(\OR)}$ and $\|\cA_0^{-1}(k)\|_{\cH^*(R+1)\to \cH(R+1)}$.
\epf

\bre\label{rem:weighted}
\mythmname{The rationale behind writing the Helmholtz equation as \eqref{eq:Helmholtz} and working in the weighted norms \eqref{eq:weighted_norms}}
When using these norms, 

(i) the solution operator has the nice property that its $k$-dependence is the same regardless of the spaces it is considered on (by Lemma \ref{lem:L2H1} and Theorem \ref{thm:H2}), and 

(ii) with the Helmholtz equation written as \eqref{eq:Helmholtz}, the continuity constants of the resulting sesquilinear form $a_0$ \eqref{eq:EDPa} and also the perturbation $a_p$ \eqref{eq:a_p} are then independent of $k$ (Parts (i) and (ii) of Lemma \ref{lem:propa}).

In the numerical-analysis literature, one usually writes the Helmholtz equation as $\nabla \cdot(\MA \nabla u) +k^2 n u=-f$ and uses the weighted $H^1$ norm
$\sqrt{\|\nabla v\|^2_{L^2} + k^2\|v\|^2_{L^2}}$. In this norm, the $k$-dependence of the solution operator as a map from $L^2\to L^2$ and as a map from $L^2 \to H^1$ is different, i.e., we lose the desirable property (i) above. One could still work in the weighted norms \eqref{eq:weighted_norms}, and write the Helmholtz equation as $\nabla \cdot(\MA \nabla u) +k^2 n u=-f$, however, the continuity constant of the sesquilinear form $a$ is then proportional to $k$, i.e., we lose the desirable property (ii) above.
\ere

\section{Proof of Theorem \ref{thm:2}}\label{sec:proof2}

\emph{Proof of Part (i).} By definition 
\beq\label{eq:Neumann1}
\cA(k,\newy)= \cA_0(k)+ \cA_p(k,\newy) = \Big(\Id + \cA_p(k,\newy) \cA_0^{-1}(k)\Big)\cA_0(k)
\eeq
so that, at least formally, 
\beq\label{eq:formal}
\cA^{-1}(k,\newy) = \cA_0^{-1}(k) \Big(\Id + \cA_p(k,\newy) \cA_0^{-1}(k)\Big)^{-1}.
\eeq
Since $\cA_0:\cH\to \cH^*$ is invertible by Theorem \ref{thm:Fred}, to
show that $\cA:\cH\to \cH^*$ is invertible, it is sufficient to show
that $\Id + \cA_p(k,\newy) \cA_0^{-1}(k)$ is invertible on $\cH^*$, and by Neumann series it
is sufficient to show that $\|\cA_p(k,\newy) \cA^{-1}_0(k)\|_{\cH^*\to \cH^*}<1$.

By Theorem \ref{thm:Fred} and Corollary \ref{cor:operators}, $\cA_p(k,\newy) \cA^{-1}_0(k): \cH^*\to \cH^*$ is bounded, and 
by the second bound in \eqref{eq:normboundA_p}
\beqs
\N{\cA_p(k,\newy) \cA^{-1}_0(k)}_{\cH^*\to \cH^*}\leq\max\Big\{ \N{\MA_p(\newy)}_{L^\infty(\Omega_R)}, \N{n_p(\newy)}_{L^\infty(\Omega_R)}\Big\} \N{\cA^{-1}_0(k)}_{\cH^*\to \cH}.
\eeqs
The condition \eqref{eq:condition1} therefore implies that $\|\cA_p(k,\newy) \cA^{-1}_0(k)\|_{\cH\to \cH^*}\leq 1/2$; thus $\Id + \cA_p(k,\newy) \cA_0^{-1}(k)$ is invertible with
\beqs
\N{(\Id + \cA_p(k,\newy) \cA_0^{-1}(k))^{-1}}_{\cH^*\to \cH^*} \leq 2,
\eeqs
and the bound \eqref{eq:thm_bound1} on $\|\cA^{-1}(k)\|_{\cH^*\to \cH}$ follow from \eqref{eq:formal}.

The bound 
\eqref{eq:normboundA_p} and the 
assumptions that the maps $\by\mapsto \MA_p(\cdot,\by)$ and $\by\mapsto n_p(\cdot,\by)$ are holomorphic for $\newy \in Y_0$ 
imply $\cA_p(k,\newy): \cH\to\cH^*$ is holomorphic for 
$\newy \in Y_0$. 
Thus $\cA_p(k,\newy) \cA_0^{-1}(k)$ is a holomorphic family in $\newy$ of bounded operators on $\cH^*$ for 
$\newy \in Y_0$. 
The claim that $\cA^{-1}(k,\newy)$ is holomorphic as a function of $\newy$ for $\newy \in Y_{\es 1}(k)$
then follows by Neumann series, since  the Neumann series converges uniformly if \eqref{eq:condition1} holds, and a uniformly-converging sum of holomorphic functions is holomorphic.

\

\noi\emph{Proof of Part (ii).} 
By Theorem \ref{thm:H2}, the assumptions that $\Omega_-$ is $C^{1,1}$ and $\MA_0\in W^{1,\infty}(\OR,\Rea^d\times \Rea^d)$ implies that $\cA^{-1}_0(k):L^2\to H^2_k(\OR)\cap \cH$.
The proof is then essentially identical to that above once we have shown that 
\begin{align*}
&\N{\cA_p(k,\newy)\cA_0^{-1}(k)}_{L^2(\OR)\to L^2(\OR)}
\\
&\hspace{3cm}\leq \max\Big\{ C\N{\MA_p(\newy)}_{W^{1,\infty}(\Omega_R)}, \N{n_p(\newy)}_{L^\infty(\Omega_R)}\Big\}
\N{\cA_0^{-1}(k)}_{L^2(\OR)\to H^2_k(\OR)\cap \cH}.
\end{align*}
To prove this inequality, by \eqref{eq:Ajdual}, it is sufficient to
show that 
\begin{align}\nonumber
&|\langle \cA_p(k,\newy)\cA_0^{-1}(k) f, v\rangle_{L^2(\OR)\times L^2(\OR)}| 
= 
|a_p(\cA_0^{-1}f,v)| \\
&
\leq  \max\Big\{ C\N{\MA_p(\newy)}_{W^{1,\infty}(\Omega_R)}, \N{n_p(\newy)}_{L^\infty(\Omega_R)}\Big\} 
\N{\cA_0^{-1}}_{L^2(\OR)\to H^2_k(\OR)\cap \cH}
\N{f}_{L^2(\OR)}\N{v}_{L^2(\OR)}
\label{eq:final}
\end{align}
for all $v$ in a dense subset of $L^2(\OR)$.
When $\MA_p\in W^{1,\infty}(\OR,\Rea^d\times \Rea^d) = C^{0,1}(\OR,\Rea^d\times \Rea^d)$ (by, e.g., \cite[\S4.2.3, Theorem 5]{EvGa:92}), $u\in H^2(\OR)$, and $v\in C^{\infty}_{\rm \comp}(\OR)$, by using integration by parts/Green's identity (see \cite[Lemma 4.1]{Mc:00}) in the definition of $\MA_p$ \eqref{eq:a_p}, we have
\beqs
a_p(u,v)= -\int_{\domain_R} 
\Big(k^{-2}\overline{v} \nabla\cdot (\MA_p\nabla u)
 + n_p u \vb\Big),
\eeqs
where we have used that $\supp\, v \Subset \OR$ so that there are no integrals on $\GR$ or $\GammaD$. 
Therefore, given $k_0>0$, there exists $C>0$ such that, for all $f\in L^2(\OR)$ and $k\geq k_0$
for all $v\in C^{\infty}_{\rm \comp}(\OR)$, 
\begin{align*}
|a_p(\cA_0^{-1}(k)f,v)|&\leq \max\Big\{ C\N{\MA_p(\newy)}_{W^{1,\infty}(\Omega_R)}, \N{n_p(\newy)}_{L^\infty(\Omega_R)}\Big\} \N{\cA_0^{-1}(k)f}_{H^2_k(\OR)} \N{v}_{L^2(\OR)},
\end{align*}
which implies \eqref{eq:final} and the proof is complete.

\section{Proofs of Theorem \ref{thm:1} and Corollary \ref{cor:1}}\label{sec:simpleproofs}

\bpf[Proof of Theorem \ref{thm:1}]
Solving the PDEs in \eqref{eq:PDE} and imposing the boundary condition at $x=0$ and the outgoing condition, we obtain that
\beqs
u(x) = A \sin\left( kx\sqrt{1/2- \parao}\right), \quad 0\leq x\leq 1, \qquad u(x) = B \re^{\ri k x}, \quad x>1.
\eeqs
For this to be a weak solution of $Pu=0$, 
we require that both $u$ and $\partial_x$ be continuous at $x=1$ 
(see, e.g., \cite[Lemma 4.19]{Mc:00}),
and this leads to the condition that 
\beq\label{eq:1}
\tan \left( k\sqrt{\half - \parao}\right) =- \ri \sqrt{\half - \parao}.
\eeq
Let $\omega(y,k) \in \Com$ (with $y\in \Com$) be such that 
\beqs
\sqrt{\half - y} = \frac{1}{\sqrt{2}} - k^{-1} \omega(y,k);
\eeqs
observe that, as $k\to \infty$, $y=\mathcal{O}(k^{-1})$ iff $\omega = \mathcal{O}(1)$.
The equation \eqref{eq:1} then becomes
\beqs
\tan\left(\frac{k}{\sqrt{2}}- \omega(y,k)\right) = -\ri  \left(\frac{1}{\sqrt{2}} - k^{-1} \omega(y,k)\right),
\eeqs
and simplifies further to 
\beq\label{eq:2}
\tan \big(\omega(y,k)\big) = \ri  \left(\frac{1}{\sqrt{2}} - k^{-1} \omega(y,k)\right),
\eeq
with $k$ as in \eqref{eq:k2}.
When $k^{-1}=0$, \eqref{eq:2} has a solution $\omega=\omega^*\approx
0.88 \ri(1+ \parat)^{-1}$. By the inverse function theorem, for $k$
sufficiently large  (i.e., for $m$ in \eqref{eq:k2} sufficiently large), \eqref{eq:2}
has a solution $\omega(y,k)$ contained in $k$-independent neighbourhood of
$\omega^*$. Moving back from the $\omega(y,k)$ variable to the $\parao$
variable, we see that there exist $C_1, C_2>0$, independent of $k$,
such that \eqref{eq:PDE} has a non-zero solution for $y$ satisfying
\eqref{eq:strip2} and $k$ sufficiently large. 
\epf

\

\bpf[Proof of Corollary \ref{cor:1}]
We show that 

(i) $\cA^{-1}(k,y): L^2([0,2])\to L^2([0,2])$ is meromorphic for $y\in \Com$, and 

(ii) if
there exists a non-zero $u\in \cH$ such that $\cA(k,y_0)u=0$, then the map $y\mapsto \cA^{-1}(k,y)$ has a pole at $y=y_0$. 

(As noted in \S\ref{sec:set2}, the proof of (i) is essentially contained in \cite[Lemma 1.12]{GaMaSp:21}, but since the proof is relatively short, we include it for completeness here.)

After having shown (i) and (ii), the result of Corollary \ref{cor:1} follows by noting that Theorem \ref{thm:2} shows that if $k$ satisfies \eqref{eq:k2} then there exists $u\in \cH$ such that, when $y=y_0$, $a(u,v)=0$ for all $v\in \cH$, and thus $\cA(k,y_0) u=0$ by \eqref{eq:Ajdual}.

For (i), by \eqref{eq:Neumann1},
it is sufficient to show that 
\beq\label{eq:toshow}
\big(\Id +\cA_p(k,y) \cA_0^{-1}(k) \big)^{-1}: L^2([0,2])\to 
L^2([0,2])
\text{ is meromorphic for } y\in \Com.
\eeq
Since $\cA_p(k,y)$ is holomorphic for $y\in \Com$,  
 the analytic Fredholm
  theorem (see, e.g., \cite[Theorem VI:14]{ReeSim72}, \cite[Theorem C.8]{DyZw:19}) implies that 
\eqref{eq:toshow} holds if 

(a) $\cA_p(k,y) \cA_0^{-1}(k): L^2([0,2])\to L^2([0,2])$ is compact (and thus $I+\cA_p(k,y) \cA_0^{-1}(k)$ is Fredholm) and 

(b) there exists $y_0\in \Com$ such that $(\Id +\cA_p(k,y) \cA_0^{-1}(k) )^{-1}$ exists.

The condition in (b) holds with $y_0=0$, since $\cA_p(k,0)=0$.
  The condition in (a) follows since $\cA_p(k, y): L^2([0,2])\to L^2([0,2])$, $\cA_0^{-1}(k): L^2([0,2])\to \cH$, and the injection $\cH\to L^2([0,2])$ is compact by
the Rellich--Kondrachov theorem; see, e.g., \cite[Theorem 3.27]{Mc:00}. 
  
For (ii), by \eqref{eq:formal}, it is sufficient to show that $(\Id + \cA_p(k,y_0) \cA_0^{-1}(k))$ has a non-empty nullspace. By assumption $\cA(k,y_0)u=0$ for $u\neq 0$; thus $(\Id +\cA_p(k,y_0) \cA_0^{-1}(k) )\widetilde{u} =0$ by \eqref{eq:Neumann1}, where $\widetilde{u} := \cA_0(k) u \neq 0$ (since $\cA_0(k)$ is invertible).  
\epf

\section{How Lemma \ref{lem:meromorphy} and Theorem \ref{thm:GMS} follow from the results of \cite{GaMaSp:21}}\label{sec:GMS}

The idea behind the results of \cite{GaMaSp:21} is the following. 

Recall that \emph{resonances} are  poles of the meromorphic continuation of the Helmholtz solution operator as a function of $k$ from $\Im k\geq 0$ to $\Im k<0$; see \cite[Theorem 4.4 and Definition 4.6]{DyZw:19}.
There has been a large amount of research showing that existence
of quasimodes (in
the sense of Definition \ref{def:quasimodesh}) with
super-algebraically-small quality implies existence of resonances super-algebraically
close to $\Im k=0$, and vice versa; see \cite{TaZw:98, St:99, St:00}
(following \cite{StVo:95, StVo:96}) and \cite[Theorem 7.6]{DyZw:19}. 

The paper \cite{GaMaSp:21} investigates eigenvalues of the truncated Helmholtz solution operator for the exterior Dirichlet problem, where $\mu\in \mathbb{C}$ is an eigenvalue of the truncated exterior Dirichlet problem at wavenumber $k>0$, with corresponding eigenfunction $u\in H^1(\OR)$ with $\gamma u=0$ on $\partial \Omega_-$, if 
\beqs
(-k^{-2}\Delta -1)u = \mu u \,\,\text{ on }\OR \quad\tand \quad\partial_{\nu} u = \DtN \gamma u \,\,\text{ on }\GR
\eeqs
(where $\DtN$ is defined in \S\ref{sec:src}); see \cite[Definition 1.1]{GaMaSp:21}.

These eigenvalues are therefore poles in $y$ of the solution operator of the problem 
\beq\label{eq:GMS2}
(-k^{-2} \Delta -1 - y)u = f  \,\,\text{ on }\OR, \quad \gamma u =0 \,\,\ton \partial \Omega_-, \,\,\tand\,\,\partial_{\nu} u = \DtN \gamma u \,\,\text{ on }\GR.
\eeq
In the notation of \S\ref{sec:set1} with 
$\MA\equiv \MI$, $n_0\equiv 1$, and $n_p(\bx,y):= y \ind_{\OR}(\bx)$, these are then poles of $y\mapsto \cA^{-1}(k,y)$.

The paper \cite{GaMaSp:21} repeats the ``quasimode-to-resonances'' arguments of \cite{TaZw:98, St:99, St:00} for the solution operator of \eqref{eq:GMS2} thought of as a function of $y$, rather than $k$; see the overview discussion in \cite[\S1.5]{GaMaSp:21}. In particular, meromorphy of the solution operator of \eqref{eq:GMS2} as a function of $y$, i.e., Lemma \ref{lem:meromorphy}, in proved in \cite[Lemma 1.12]{GaMaSp:21}, with the ``quasimodes to eigenvalues'' (i.e., ``quasimodes to poles in $y$'') result of Theorem \ref{thm:GMS} proved in \cite[Theorems 1.5 and 2.2]{GaMaSp:21}.

While \cite{GaMaSp:21} is concerned with the eigenvalues of the constant-coefficient Helmholtz equation outside a Dirichlet obstacle $\Omega_-$, we expect the arguments to carry over to \eqref{eq:Helmholtz} with $\Omega_-=\emptyset$ and smooth $\MA_0$ and $n_0$ supporting quasimodes with superalgebraically-small quality. Indeed, the key result to prove is the bound in  \cite[Lemma 1.14]{GaMaSp:21} on the solution operator away from the poles in $y$, and the proof of this in the boundaryless case should be easier than the proof with boundary in \cite[\S3]{GaMaSp:21}.

\section*{Acknowledgements}

Both authors thank the anonymous referees for their constructive comments on a previous version of the paper.
EAS thanks 
Christoph Schwab (ETH Z\"urich) for useful discussions on parametric holomorphy in UQ.
EAS  was supported by EPSRC grant EP/1025995/1.
JW was partially supported by Simons Foundation grant 631302, NSF grant DMS--2054424, and a Simons Fellowship.

\footnotesize{
\bibliographystyle{plain}
\bibliography{biblio_combined_sncwadditions}
}

\end{document}